\documentclass{conm-p-l}
\newtheorem{theorem}{Theorem}[section]
\newtheorem{lemma}[theorem]{Lemma}

\theoremstyle{definition}

\theoremstyle{remark}
\newtheorem{remark}[theorem]{Remark}

\numberwithin{equation}{section}
\usepackage{epsfig}
\newtheorem{proposition}[theorem]{Proposition}

\newcommand{\cal}{\mathcal}

\newcommand{\CC}{{\mathbb{C}}}
\newcommand{\HH}{{\mathbb{H}}}
\newcommand{\FF}{{\mathbb{F}}}
\newcommand{\PP}{{\mathbb{P}}}
\newcommand{\QQ}{{\mathbb{Q}}}
\newcommand{\RR}{{\mathbb{R}}}
\newcommand{\ZZ}{{\mathbb{Z}}}

\renewcommand{\phi}{\varphi}
\begin{document}
\title{Degenerations of $(1,3)$ abelian surfaces and Kummer surfaces}
\author{K. Hulek} 

\address{Institut f\"ur Mathematik, Universit\"at Hannover, Postfach 6009,
D~30060 Hannover, GERMANY}

\email{hulek@math.uni-hannover.de}

\author{I. Nieto}

\address{Cimat, A.C., Callejon de Jalisco Sn, col. Mineral de Valenciana,
36000~Guanajuato, Gto., MEXICO}

\curraddr{Instituto de Fisica y Matematicas, Universidad Michoacan de San
Nicolas Hidalgo, Morelia, Michoacan, MEXICO}

\email{nieto@fractal.cimat.mx}

\thanks{The second author was supported by the UMSNH research project
``Modelos proyectivos asociados a Variedades degeneradas de Kummer'' and by
EPSRC grant number GR/L27534.}

\author{G.K. Sankaran} 

\address{Department of Mathematical Sciences, University of Bath, Claverton
Down, Bath BA2~7AY, ENGLAND}

\email{gks@maths.bath.ac.uk} 

\thanks{The first and third authors were partially supported by EU HCM
network AGE ERBCHRX-CT94-0557.}

\date{26th February 1999}

\begin{abstract}
We study degenerations of Kummer surfaces associated to certain
divisors in Nieto's quintic threefold and show how they arise from
boundary components of a suitable toroidal compactification of the
corresponding Siegel modular threefold.
\end{abstract}

\maketitle

The aim of this paper is to study the degenerations of an interesting class
of Kummer surfaces in $\PP^3$ in terms of degenerations of the
corresponding abelian surfaces. The quintic threefold in $\PP^4$ given by
$$
N=\left\{ \sum\limits^5_{i=0} u_i=\sum\limits^5_{i=0}\frac 1{u_i}=0\right\}
\subset
\PP^5.
$$
is the closure of the locus parametrizing $H_{22}$-invariant quartics
with $16$ skew lines. The smooth surfaces of this type are parametrized by
a non-empty open set $N^s$ of $N$. These surfaces are Kummer surfaces
associated to abelian surfaces with a $(1,3)$--polarization. The action of
the Heisenberg group on the Kummer surface corresponds to a level-$2$
structure on the abelian surface. If $A$ is an abelian surface and ${\cal
L}$ is a symmetric line bundle representing this polarization then ${\cal
L}^{\otimes 2}$ is the unique totally symmetric line bundle representing
the $(2,6)$--polarization. Under the involution $\iota:x\mapsto -x$ the
space $H^0({\cal L}^{\otimes 2})$ decomposes into eigenspaces $H^0({\cal
L}^{\otimes 2})^+$ and $H^0({\cal L}^{\otimes 2})^-$ of dimensions~$8$
and~$4$ respectively. The linear system given by $H^0({\cal L}^{\otimes
2})^-$ defines a rational map $A--\rightarrow \PP^3$ which induces a map
from the smooth Kummer surface $\widetilde{\operatorname{Km}}(A)$ to
$\PP^3$. In fact for general $A$ this defines an embedding of
$\widetilde{\operatorname{Km}}(A)$ into $\PP^3$. With respect to a suitable
basis of the linear system which defines the map to $\PP^3$ the image
surface is $H_{22}$-invariant.

Let ${\cal A}_{1,3}(2)$ be the moduli space of $(1,3)$--polarized abelian
surfaces with a level-$2$ structure. Then the map which associates to an
abelian surface its Kummer surface defines a $2:1$ map ${\cal
A}_{1,3}(2)\rightarrow N$ (see also the paper \cite{HNS}). It was already
shown by one of the authors \cite{Ni} that the complement of $N^s$ in $N$
consists of two sets of $15$ planes called the V- and S-planes respectively
(see also \cite{BN}). In fact this author studied in \cite{Ni} the
configuration of V- and S-planes and computed the structure of the singular
Kummer surface which belongs to a general point on these planes showing in
particular that the generic point on a V-plane corresponds to a quotient of
a product of elliptic curves by a $\ZZ/2 \times \ZZ/2$ action. In
\cite{HNS} we proved this result by a different method and showed that the
Kummer surfaces which are parametrized by the V-planes are exactly the
bielliptic surfaces in the sense of \cite{HW}. In this paper we consider
the quartic surfaces which are parametrized by the $15$ S-planes and we
will show that they correspond to Kummer surfaces of degenerate abelian
surfaces. More precisely we construct an explicit degeneration of abelian
surfaces which gives rise to these singular Kummer surfaces.

\section{Theta functions}
In this section we will give an explicit description of a basis of the
space $H^0({\cal L}^{\otimes 2})^-$ which defines the map from $A$ to
$\PP^3$, factoring through the Kummer surface, in terms of theta functions.
Our standard reference for theta functions is Igusa's book \cite {I}. We
shall denote points of the Siegel upper half plane $\HH_2$ by $\tau=\left(
\begin{array}{cc}
\tau_1 & \tau_2\\
\tau_2 & \tau_3
\end{array}
\right),$ and $z=(z_1, z_2)$ will denote the coordinates on $\CC^2$. For
every pair $(m', m'')\in\RR^2\times\RR^2$ we define the theta function
$$
\Theta_{m' m''}(\tau, z)=\sum\limits_{q\in\ZZ^2} e^{2\pi i[\frac 12
(q+m')\tau^t(q+m')+(q+m')^t(z+m'')]}.
$$
Given a point $\tau=\left(
\begin{array}{cc}
\tau_1 & \tau_2\\
\tau_2 & \tau_3
\end{array}
\right)\in \HH_2$ we associate to it a period matrix
$$
\Omega_{\tau}=\left(
\begin{array}{cccc}
2\tau_1 & 2\tau_2 & 2 & 0\\
2\tau_2 & 2\tau_3 & 0 & 6
\end{array}
\right)
$$
and the lattice
$$
L_{\tau}=\ZZ^4\Omega_{\tau}=\ZZ e_1+\ZZ e_2 + \ZZ e_3 +\ZZ e_4
$$
generated by the columns $e_i$ of the period matrix $\Omega_{\tau}$. The
abelian surface
$$
A_{\tau}=\CC^2/L_{\tau}
$$
has a $(1,3)$--(and hence also a $(2,6)$--)polarization. Normally $0\in
A_{\tau}$ is chosen as the origin and the involution given by taking the
inverse is $\iota:x\mapsto -x$. For reasons which will become apparent
later we shall want to define the origin of $A_{\tau}$ as the image of the
point
$$
\omega=\frac 12(1,1)\left(
\begin{array}{cc}
\tau_1 & \tau_2\\
\tau_2 & \tau_3
\end{array}
\right)=\frac 12 (\tau_1+\tau_2, \tau_2+\tau_3).
$$
Note that with respect to $0$ this is a $4$-torsion point. Then the
involution with respect to $\omega$ is given by
$$
\iota_{\omega}(z)=-z+2\omega.
$$
Finally we set
$$
\tau'=\left(
\begin{array}{cc}
\tau_1/2  & \tau_2/6\\
\tau_2/6  &  \tau_3/18
\end{array}
\right),\quad
z'=(z_1/2,z_2/6).
$$
The main objects of this section are the functions
$$
{\widehat{\Theta}}_{\alpha \beta}(\tau,z): =
\Theta_{00\frac{\alpha}{2}\frac{\beta}{6}}(\tau',z'-{\omega}'),
\quad \alpha=0,1;\
\beta=0,\ldots,5.
$$
Note that ${\widehat{\Theta}}_{\alpha+2, \beta}={\widehat{\Theta}}_{\alpha
\beta}$ and $\widehat{\Theta}_{\alpha, \beta+6}=\widehat{\Theta}_{\alpha
\beta}$
so that we can read the indices cyclically.

\begin{lemma}\label{lem11}
{\rm(i)} The functions ${\widehat{\Theta}}_{\alpha\beta}$ are
all sections of the same line bundle ${\cal L}_{\tau}$ on $A_{\tau}$.\\
{\rm(ii)} ${\cal L}_{\tau}$ represents a polarization of type $(2,6)$.
\end{lemma}

\begin{proof}
  (i)\ We must prove that the automorphy factor of the functions
  ${\widehat{\Theta}}_{\alpha \beta}$ with respect to $z\mapsto z+e_i$ does
  not depend on $(\alpha, \beta)$. This follows immediately from the
  formulae $(\Theta 1)-(\Theta 5)$ of \cite[pp. 49, 50]{I}.\\
\noindent
(ii)\ Since the type of a polarization is constant in families it is enough
to prove the statement for $\tau_2=0$ where
$$
A_{\tau}=E({\tau_1}) \times E({\tau_3})
$$
with
$$
E(\tau_1)=\CC/(\ZZ 2\tau_1+\ZZ 2),\ E(\tau_3)=\CC/(\ZZ 2\tau_3+\ZZ 6).
$$
In this case
$$
{\widehat{\Theta}}_{\alpha
\beta}(\tau,z)=\vartheta_{0\frac{\alpha}{2}}(\tau_1/2, z_1/2-{\omega_1}/2)\
\vartheta_{0\frac{\beta}{6}}(\tau_3/18, z_2/6-{\omega_2}/6)
$$
where we use $\vartheta$ to denote theta functions in one variable. We
claim that the degree on $E(\tau_1)$ is $2$ and that the degree on
$E(\tau_3)$ is $6$. Indeed the first claim follows since
$$
\begin{array}{rcl}
\vartheta_{0\frac{\alpha}{2}}(\tau_1/2, z_1/2-{\omega_1}/2)=0 &
\Leftrightarrow &
z_1/2 \in (\ZZ \tau_1/2 +\ZZ)-\alpha/2+{\omega_1}/2\\
&\Leftrightarrow & z_1\in (\ZZ \tau_1+\ZZ 2)-\alpha+{\omega_1}.
\end{array}
$$
This means that $\vartheta_{0\frac{\alpha}{2}}(\tau_1/2,
z_1/2-{\omega_1}/2)$ has two zeroes on $E(\tau_1)$. The other claim follows
in exactly the same way.  \hfill
\end{proof}

We shall denote the sections of ${\cal L}_{\tau}$ defined by
${\widehat{\Theta}}_{\alpha\beta}(\tau, z)$ by ${\widehat{s}}_{\alpha
\beta}$. By general theory the twelve sections
${\widehat{s}}_{\alpha\beta};\quad \alpha=0,1,\ \beta=0,\ldots,5$ form a
{\em basis} of $H^0({\cal L}_{\tau})$. (Cf. \cite[p.75]{I} for an analogous
statement.)

We now want to describe the symmetry properties of the line bundle ${\cal
L}_{\tau}$ and the sections ${\widehat{s}}_{\alpha\beta}$. The kernel of
the map
$$
\begin{array}{rcl}
\lambda:\ A_{\tau}  &  \rightarrow  &  \mbox{Pic}^0A_{\tau}\\
x & \mapsto & t_x^*{\cal L}_{\tau} \otimes{\cal L}_{\tau}^{-1}
\end{array}
$$
where $t_x$ is translation by $x$ is equal to
$$
\operatorname{ker } \lambda =(\ZZ \frac{e_1}{2}
+\ZZ\frac{e_2}{6}+\ZZ\frac{e_3}{2}+\ZZ\frac{e_4}{6}) L_{\tau}\cong(\ZZ/2)^2
\times(\ZZ/6)^2.
$$
\noindent
We set $\rho_6:=e^{2\pi i/6}$.

\begin{proposition}\label{prop12}
{\rm{(i)}} The group $\operatorname{ker} \lambda$ acts on the sections
${\widehat{s}}_{\alpha \beta}$ as follows
$$
\begin{array}{rlcllclcl}
e_1 /2 & : & {\widehat{s}}_{\alpha\beta} & \mapsto &
(-1)^{\alpha}\ {\widehat{s}}_{\alpha\beta}\ , & e_2 /6 & :
{\widehat{s}}_{\alpha\beta} & \mapsto & \rho_6^{-\beta}\
{\widehat{s}}_{\alpha\beta}\\
e_3 /2 & : & {\widehat{s}}_{\alpha\beta} & \mapsto &
{\widehat{s}}_{\alpha+1, \beta}\ , & e_4 /6 & :
{\widehat{s}}_{\alpha\beta} & \mapsto & {\widehat{s}}_{\alpha, \beta+1}.
\end{array}
$$
{\rm{(ii)}} The involution $\iota_{\omega}$ acts on the sections
${\widehat{s}}_{\alpha\beta}$ by
$$
\iota_{\omega}:{\widehat{s}}_{\alpha\beta}\mapsto {\widehat{s}}_{-\alpha,
-\beta}.
$$
\end{proposition}

\begin{proof}
  (i)\ We shall prove this for $e_1/2$, the other cases being similar.
  Again using \cite[pp. 49, 50]{I} we find

$$
\begin{array}{rcl}
{\widehat{\Theta}}_{\alpha\beta}(\tau, z+\frac{e_1}{2}) & =
&  \Theta_{00\frac{\alpha}{2}\frac{\beta}{6}}(\tau',
z'-{\omega}'+(\tau_1/2, \tau_2/6))\\[2mm]
& {\begin{array}{c}
(\Theta 3)\\
{\sim}
\end{array}}
&  e^{2\pi i (-(1,0)({\begin{array}{c}\alpha/2\\ \beta/6
\end{array})})} \Theta_{10\frac{\alpha}{2}\frac{\beta}{6}}(\tau',
z'-{\omega}')\\[2mm]
& {\begin{array}{c}
(\Theta 1)\\
=
\end{array}} &  (-1)^{\alpha} \Theta_{00\frac{\alpha}{2}\frac{\beta}{6}}(\tau',
z'-{\omega}')\\ [2mm]
& =  &  (-1)^{\alpha}{\widehat{\Theta}}_{\alpha\beta}(\tau, z).
\end{array}
$$
\noindent
Here $\sim$ denotes equality up to a nowhere vanishing function which is
independent of $\alpha$ and $\beta$.\\
\noindent (ii) Here we have that
$$
\begin{array}{rcl}
{\widehat{\Theta}}_{\alpha\beta}(\tau, -z+2\omega) & = &
\Theta_{00\frac{\alpha}{2}\frac{\beta}{6}}(\tau', -z'+{\omega}')\\
& {\begin{array}{c}
(\Theta 1)\\
=
\end{array}} &  \Theta_{00-\frac{\alpha}{2}-\frac{\beta}{6}}(\tau',
z'-{\omega}')\\ & = & {\widehat{\Theta}}_{-\alpha -\beta}(\tau, z)
\end{array}
$$
where indices are to be read cyclically.\hfill
\end{proof}

\begin{remark}\label{rem13}
  {\rm{(i)}} Part (i) of the above proposition gives an explicit
  description of the lifting of the group $(\ZZ /2)^2 \times (\ZZ /6)^2$ to
  the Heisenberg group $H_{26}$.\\
  {\rm{(ii)}} Note that part (ii) of the above proposition is true for any
  choice of the point $\omega$ and hence in particular also for the
  involution $\iota$ itself.
\end{remark}

We can now describe a basis of the eigenspaces
$H^0({\cal L}_\tau)^+$ and
$H^0({\cal L}_\tau)^-$ as follows:
$$
\begin{array}{ll}
{\widehat{u}}_{\alpha\beta}={\widehat{s}}_{\alpha\beta}+{\widehat{s}}
_{-\alpha,-\beta}\in H^0({\cal L}_\tau)^+;  &
\alpha \in \{0,1\}, \beta \in \{0,1,2,3\}\\
{\widehat{t}}_{\alpha\beta}={\widehat{s}}_{\alpha\beta}-{\widehat{s}}
_{-\alpha,-\beta}\in H^0({\cal L}_\tau)^-;  &
(\alpha,\beta)= (0,1),(0,2),(1,1),(1,2).
\end{array}
$$
\noindent
For our purposes it is, however, better to work with a different basis of
$H^0({\cal L}_\tau)^-$.
$$
\begin{array}{ccrclclcl}
{\widehat{g}}_0 & := & {\widehat{t}}_{01} & + &{\widehat{t}}_{11}& - &
{\widehat{t}}_{02} & - & {\widehat{t}}_{12}\\
{\widehat{g}}_1 & := & -{\widehat{t}}_{01} & - &{\widehat{t}}_{11}& - &
{\widehat{t}}_{02} & - & {\widehat{t}}_{12}\\
{\widehat{g}}_2 & := & {\widehat{t}}_{01} & - &{\widehat{t}}_{11}& - &
{\widehat{t}}_{02} & + & {\widehat{t}}_{12}\\
{\widehat{g}}_3 & := & -{\widehat{t}}_{01} & + &{\widehat{t}}_{11}& - &
{\widehat{t}}_{02} & + & {\widehat{t}}_{12}.
\end{array}
$$
Recall the Heisenberg group $H_{22}$ from \cite {BN}. The group $H_{22}$ has
order~$32$ and
$$
H_{22}/\mbox{ centre } \cong (\ZZ/2)^4
$$
is the group generated by the elements
$$
\begin{array}{cclcl}
\sigma_1 & : & (z_0: z_1: z_2: z_3) & \mapsto & (z_2: z_3: z_0: z_1)\\
\sigma_2 & : & (z_0: z_1: z_2: z_3) & \mapsto & (z_1: z_0: z_3: z_2)\\
\tau_1   & : & (z_0: z_1: z_2: z_3) & \mapsto & (z_0: z_1: -z_2:-z_3)\\
\tau_2   & : & (z_0: z_1: z_2: z_3) & \mapsto & (z_0:-z_1: z_2: -z_3).
\end{array}
$$
The group ${}_2A_{\tau}$ of $2$-torsion points of $A_{\tau}$ is
contained in $\operatorname{ker} \lambda$. Here we identify ${}_2A_{\tau}$
with translations of $A_{\tau}$ of order~$2$. Using the translations
$x\mapsto x+e_i/2$ as generators we obtain an identification of
${}_2A_{\tau}$ with $(\ZZ/2)^4$. A straightforward calculation using
Proposition \ref{prop12} and the definition of the basis
${\widehat{g}}_0,\ldots,{\widehat{g}}_3$ shows that
$$
\begin{array}{clclll}
e_1/2 : &(\widehat{g}_0: \widehat{g}_1: \widehat{g}_2:
\widehat{g}_3) & \mapsto &
(\widehat{g}_2: \widehat{g}_3: \widehat{g}_0: \widehat{g}_1) & = &
\sigma_1 (\widehat{g}_0:
\widehat{g}_1: \widehat{g}_2: \widehat{g}_3)\\
e_2/2  : &(\widehat{g}_0: \widehat{g}_1: \widehat{g}_2:
\widehat{g}_3) & \mapsto &
(\widehat{g}_1: \widehat{g}_0: \widehat{g}_3: \widehat{g}_2) & = & \sigma_2
(\widehat{g}_0:
\widehat{g}_1: \widehat{g}_2: \widehat{g}_3)\\
e_3/2  : &(\widehat{g}_0: \widehat{g}_1: \widehat{g}_2:
\widehat{g}_3) & \mapsto &
(\widehat{g}_0: \widehat{g}_1: -\widehat{g}_2: -\widehat{g}_3) & = & \tau_1
(\widehat{g}_0:
\widehat{g}_1: \widehat{g}_2: \widehat{g}_3)\\
e_4/2  : &(\widehat{g}_0: \widehat{g}_1: \widehat{g}_2:
\widehat{g}_3) & \mapsto &
(\widehat{g}_0: -\widehat{g}_1: \widehat{g}_2: -\widehat{g}_3) & = & \tau_2
(\widehat{g}_0:
\widehat{g}_1: \widehat{g}_2: \widehat{g}_3)
\end{array}
$$

We can, therefore, summarize our results as follows:

\begin{theorem}\label{theo13}
  {\rm(i)} The basis $\widehat{g}_0,\ldots,\widehat{g}_3 \in H^0({\cal
  L}_{\tau})^-$ defines a rational map from $A_{\tau}$ to $\PP^3$ which
  factors through $\widetilde{\operatorname{Km}}(A_{\tau})$. This map is
  equivariant with respect to the action of ${}_2A_{\tau}\cong (\ZZ/2)^4$
  on $A_{\tau}$ and of $H_{22}/\mbox{centre }\cong (\ZZ/2)^4$ on $\PP^3$.
  In particular the image is $H_{22}$-invariant.\\
  {\rm(ii)} The Kummer surface $\widetilde{\operatorname{Km}}(A_{\tau})$ is
  embedded as a smooth quartic surface if and only if $A_{\tau}$ is neither
  a product nor a bielliptic abelian surface.  If $A_{\tau}$ is bielliptic
  then $\widetilde{\operatorname{Km}}(A_{\tau})$ is mapped to a quartic
  with four nodes; if $A_{\tau}$ is a product, then
  $\widetilde{\operatorname{Km}}(A_{\tau})$ is mapped $2:1$ onto a quadric.
\end{theorem}

\begin{proof}
{\rm(i)} Follows immediately from our above calculations.\\
{\rm(ii)} This was shown in \cite{HNS}.
\hfill
\end{proof}

\begin{remark}
  ${\cal L}_{\tau}$ is the unique totally symmetric line bundle with
  respect to the involution ${\iota}_{\omega}$.
\end{remark}

\section{Degenerations}
In this section we construct degenerations of $(1,3)$--polarized abelian
surfaces which correspond to points on the S-planes. The construction of
degenerating families of abelian varieties is in general technically
complicated (see e.g. \cite{FC}, \cite{AN}). Although we cannot avoid these
technicalities entirely, we have tried to present our construction in a way
which uses only a minimum of technical steps. These, however, cannot be
avoided.

We consider the group
$$
P=\left\{
\left(
\begin{array}{c|c}
\left.{\bf 1}\right. &
\begin{array}{rr}
2 \ZZ &  6\ZZ\\
6 \ZZ & 18 \ZZ
\end{array}\\
\hline
0  &  {\bf 1}
\end{array}
\right)
\right\} \subset \mbox{Sp}(4,\ZZ).
$$
Note that this is the lattice contained in the parabolic subgroup of
$\Gamma_{1,3}(2)\cap \Gamma_{1,3}^{\mbox{\scriptsize lev}} $ which fixes
the isotropic plane $h=(0, 0, 1, 0)\wedge (0, 0, 0, 1)$. Here
$\Gamma_{1,3}(2)$ is the group which defines the moduli space of abelian
surfaces with a $(1,3)$--polarization and a level-$2$ structure, whereas
$\Gamma_{1,3}^{\mbox{\scriptsize lev}}$ belongs to the moduli space of
$(1,3)$--polarized abelian surfaces with a canonical level structure (cf.
\cite[I.1] {HKW}). There are two reasons for considering this group. One is
that we can then make use of the constructions in \cite{HKW} which from our
point of view is the most economical way to construct the degenerations
which we are interested in; the second reason is that, at least with the
known constructions of degenerations of abelian surfaces, the presence of a
canonical level structure is necessary. We could also use the method of
Alexeev and Nakamura \cite{AN} which likewise goes back to Mumford's
construction \cite{M}, and \cite{Nak}, \cite{Nam}. For the surfaces which
we are interested in it makes, however, little difference which of these
methods we choose.

The group $P$ acts on $\HH_2$ by
$$
\left(
\begin{array}{cc}
\tau_1 & \tau_2\\
\tau_2 & \tau_3
\end{array}
\right)\mapsto\left(
\begin{array}{cc}
\tau_1+2\ZZ & \tau_2+6\ZZ\\
\tau_2+6\ZZ & \tau_3+18\ZZ
\end{array}\right).
$$
The partial quotient of $\HH_2$ by $P$ is given by
$$
\begin{array}{ccl}
\HH_2 & \rightarrow & \HH_2/P\subset (\CC^*)^3\\
\left(
\begin{array}{cc}
\tau_1 & \tau_2\\
\tau_2 & \tau_3
\end{array}
\right) & \mapsto &
(e^{2\pi i \tau_1/2}, e^{2\pi i \tau_2/6},e^{2\pi i \tau_3/18})=(t_1,t_2,t_3).
\end{array}
$$
Recall that
$$
A_{\tau}=\CC^2/L_{\tau}
$$
where
$$
L_{\tau}=\ZZ
\left(
\begin{array}{l}
2\tau_1\\
2\tau_2
\end{array}
\right)+
\ZZ
\left(
\begin{array}{l}
2\tau_2\\
2\tau_3
\end{array}
\right)+
\ZZ
\left(
\begin{array}{c}
2\\
0
\end{array}
\right)+
\ZZ
\left(
\begin{array}{c}
0\\
6
\end{array}
\right)=L'_{\tau}+L''.
$$
Here $L_{\tau}'$ is spanned by the first two columns of $\Omega_{\tau}$
and $L''$ by the last two. Obviously $L''$ does not depend on $\tau$ and
$$
\CC^2/L''=(\CC^*)^2.
$$
We shall use the coordinates
$$
w_1=z_1/2,\quad w_2=z_2/6
$$
on $(\CC^*)^2$. The lattice $L'_{\tau}$ acts on the trivial torus bundle
$\HH_2/P\times(\CC^*)^2$ by
$$
(m, n): (t_1, t_2, t_3; w_1, w_2)\mapsto (t_1, t_2, t_3;\  t_1^{2m}\  t_2^{6n}
w_1, t_2^{2m}\  t_3^{6n} w_2).
$$
We have to extend the trivial bundle $\HH_2/P\times(\CC^*)^2$ to the
boundary in such a way that the action of $L'_{\tau}$ also extends. The
general theory of toroidal compactifications of moduli spaces of abelian
surfaces (the material which is relevant in our situation can be found in
\cite{HKW}) leads us to consider first the map
$$
\begin{array}{ccl}
(\CC^*)^3       & \rightarrow & \CC^3\\
(t_1, t_2, t_3) & \mapsto     & (t_1 t_2, t_2 t_3, t_2^{-1})=(T_1, T_2, T_3).
\end{array}
$$
Let
$$
B:=\overset{\circ}{(\overline{\HH_2/P})}
$$
be the interior of the closure of $\HH_2/P$ in $\CC^3$ in the
$\CC$--topology. (What we have considered here is an open part of the
partial compactification in the direction of the cusp corresponding to
$h=(0,0,1,0) \wedge (0,0,0,1)$. The surfaces $B\cap\{T_i=0\}$ are mapped to
boundary surfaces in the Igusa compactification of the moduli space ${\cal
A}_{1,3}^{\mbox{\scriptsize lev}}(2)$ of $(1,3)$--polarized abelian
surfaces with both a level-$2$ and a canonical level structure.) In terms
of the coordinates $T_i$ the action of $L'_{\tau}$ is now given by
$$
(m, n): (T_1, T_2, T_3; w_1, w_2)\mapsto (T_1, T_2, T_3;\  T_1^{2m} T_3^{2m-6n}
w_1, T_2^{6n} T_3^{6n-2m} w_2).
$$
Here we are particularly interested in degenerations which are given by
$\tau_1\rightarrow i\infty$. This is equivalent to $t_1=0$ and hence
corresponds to points on the surface $T_1=0$.

We now consider the space
$$
{\tilde {\cal P}}=\mbox{Proj }
R_{\Phi,\Sigma}\rightarrow \mbox{Spec } \CC[T_1, T_2,
T_3]\cong \CC^3
$$
which was defined in \cite[p.210]{HKW}. Let
$$
{\cal P}:={\tilde {\cal P}}|_{B}.
$$
Then ${\cal P}$ is a partial compactification of the trivial torus
bundle $\HH_2/P\times(\CC^*)^2$ over $\HH_2/P\subset B$. Moreover the
action of $L_{\tau}$ on the trivial torus bundle extends to an action on
${\cal P}$. The construction of ${\tilde {\cal P}}$ is originally due to
Mumford \cite[final example]{M}. Let
$$
{\bar A}:={\cal P}/L_{\tau}.
$$
Then we have a diagram

$$
\begin{array}{ccc}
A             &  \subset  &  {\bar A}\\
\pi\downarrow &           &  \downarrow\pi\\
\HH_2/P       &  \subset  &  B
\end{array}
$$
where $A=(\HH_2/P\times(\CC^*)^2)/L_{\tau}$ is the universal family. In
particular ${\bar A}$ extends the universal family $A$ to the boundary. The
fibres
$$
{\bar A}_{u}=\pi^{-1}(u)
$$
over ``boundary points'' $u\in B\backslash(\HH_2/P)$ are degenerate
abelian surfaces. We are interested in the fibres ${\bar A}_u$ over points
$u=(0, T_2, T_3)$ with $T_2 T_3\neq 0$. These are the corank~$1$
degenerations associated to the boundary component given by
$\tau_1\rightarrow i\infty$. Note that if $T_2T_3\neq 0$ then this gives
$t_2=T_3^{-1}$ and $t_3=T_2T_3$. In particular the point $u$ determines a
point $(\tau_2, \tau_3)\in\CC\times \HH_1$ where $\tau_2$ and $\tau_3$
are uniquely defined up to $6\ZZ$ and $18\ZZ$ respectively.

We can now formulate the main result of this section.

\begin{theorem}\label{theo21}
  Let $u=(0, T_2, T_3)\in B$. Then ${\bar A}_u$ is a degenerate abelian
  surface with the following properties:\\
  {\rm(i)} ${\bar A}_u$ is a corank $1$ degeneration. More precisely ${\bar
  A}_u$ is a chain of two elliptic ruled surfaces $A_{u,1}, A_{u,2}$ i.e.
  there exists an elliptic curve $E$ and a line bundle ${\cal M}_u\in
  \operatorname{Pic}^{0}(E)$ such that $A_{u,i}=\PP({\cal O}_E\oplus {\cal
  M}_u), \ i=1,2$.  The surfaces $A_{u,i}$ are glued with
  a glueing parameter $e$ as shown below in Figure~$1$.\\
  {\rm(ii)} The base curve $E\cong E(\tau_3)=\CC/(\ZZ 2\tau_3 + \ZZ 6)$.\\
  {\rm(iii)} The line bundle ${\cal M}_u={\cal
  O}_{E(\tau_3)}(6[\tau_2]-6[0])$ where $[\tau_2],[0]$
  are the points of $E(\tau_3)$ given by $\tau_2$ and~$0$.\\
  {\rm(iv)} The glueing parameter $e=[2\tau_2] \in E(\tau_3).$

\unitlength1cm
\begin{figure}[htb]
 \begin{picture}(13.5,8.5)
  \put(-0.7,0){\includegraphics{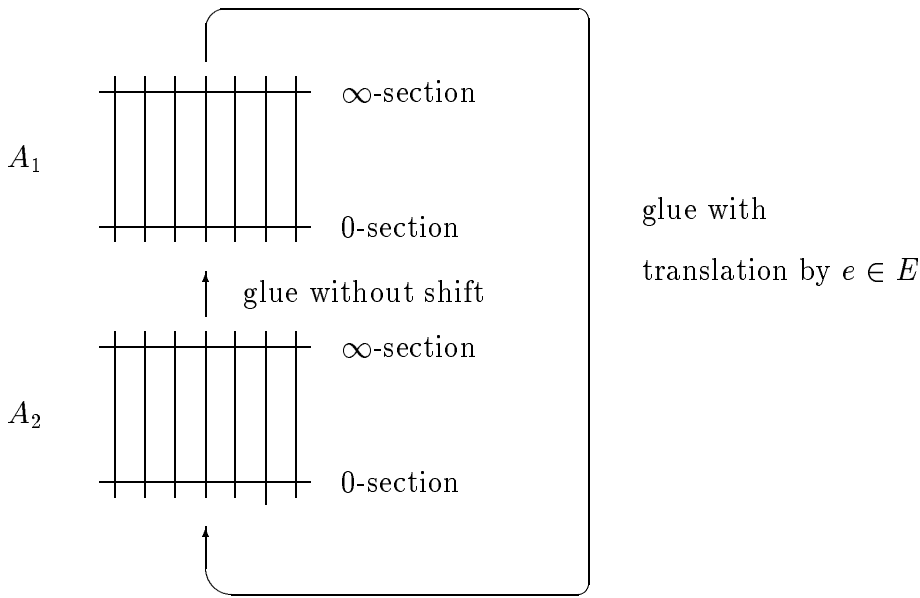}}
 \end{picture}
 \caption{Glueing of the surface $\bar A_u$}
\end{figure}

\end{theorem}
\begin{proof}
  We can derive this from \cite[part II]{HKW}. There the quotient ${\hat
  A}={\cal P}/{\hat L}$ was considered where ${\hat L}\cong \ZZ^2$ acts on
  the trivial torus bundle $\HH_2/P\times(\CC^*)^2$ by
$$
(m, n): (T_1, T_2, T_3; w_1, w_2) \mapsto (T_1, T_2, T_3;  T_1^{m} T_3^{m-n}
w_1, T_2^{n} T_3^{n-m} w_2).
$$
Hence $L'_{\tau}$ is a subgroup of ${\hat L}$ with ${\hat
L}/L_{\tau}'\cong(\ZZ/2)\times(\ZZ/6)$. This means that we can use the
description of ${\hat A}$ given in \cite[part II]{HKW} to give a
description of ${\bar A}$. In the terminology of \cite{HKW} the group
$L_{\tau}'=<s^{-2}, r^{-6}>$. The statements (i) and (ii) now follow
exactly as in the proof of \cite [Theorem II.3.10]{HKW}. In particular, the
fact that $s^2\in L_{\tau}'$, but $s \not\in L_{\tau}'$ implies that ${\bar
A}_u$ has two irreducible components. The statement about the base curve
$E$ follows from diagram \cite[II.3.13]{HKW}. Statements (iii) and (iv) are
an immediate consequence of the proof of \cite[Proposition (II.3.20)]{HKW}.
\end{proof}

Our next task is to study the involutions $\iota$ and $\iota_{\omega}$ on
$A$ and their extensions to ${\bar A}$. If we choose $0\in
A_{\tau}=\CC^2/L_{\tau}$ as the origin, then this defines a section of $A$
which extends to a section of ${\bar A}$. Moreover, the involution
$\iota:z\mapsto -z$ defines an involution of $A$ which extends to ${\bar
A}$ (this is the involution given by \cite [Lemma (II.2.9)(ii)]{HKW}. But
we said in section 1 that we wanted to choose $\omega=[(\tau_1+\tau_2)/2,
(\tau_2 + \tau_3)/2]$ as the origin. This point is a $4$-torsion point of
$A_{\tau}$ if we choose~$0$ as the origin. This choice of origin will be
necessary for what follows, but at this point it has the disadvantage that
it only defines a multisection of $A$, not a section. Nevertheless this
multisection extends to ${\bar A}$. We also claim that the involution
$\iota_{\omega}(z)=-z+2\omega$ extends to ${\bar A}$. Since $z\mapsto -z$
is defined on ${\bar A}$ it is enough to show that the translation
$z\mapsto z+2\omega$ is defined on $A$ and extends to ${\bar A}$. This is
easy to see, since $z\mapsto z+2\omega$ in terms of the coordinates $w_1,
w_2$ is given by
$$
(w_1, w_2)\mapsto(t_1 t_2^3 w_1, t_2 t_3^3 w_2)=(T_1 T_3^{-2} w_1, T_2^3 T_3^2
w_2).
$$
This is the element $s^{-1} r^{-3} \in {\hat L}$ and hence acts on
${\cal P}$ and on ${\bar A}={\cal P}/L'_{\tau}$. Since $s^{-2} r^{-6}\in
L'_{\tau}$, this is an involution. In particular $\iota_{\omega}$ defines
an involution on the fibres ${\bar A}_u$ of ${\bar A}$. Recall that for
$u=(0, T_2, T_3)$ with $T_2 T_3\neq 0$ the surface ${\bar A}_u$ has two
irreducible components $A_{u,i},\ i=1,2$ and that the singular locus of
${\bar A}_u$ consists of two disjoint elliptic curves $E_1$ and $E_2$ with
$E_i \cong E(\tau_3)=\CC/(\ZZ 2\tau_3+\ZZ 6)$.

\begin{proposition}\label{prop22}
  The involution $\iota_{\omega}$ interchanges the two components $A_{u,1}$
  and $A_{u,2}$ of ${\bar A}_u$ and induces an involution on each of the
  two curves $E_1$ and $E_2$ with four fixed points on each of these curves.
\end{proposition}

\begin{proof}
  The involution $\iota$ fixes each of the surfaces $A_{u,1}$ and $A_{u,2}$
  and interchanges $E_1$ and $E_2$. Addition by the $2$-torsion point
  $2\omega$ also interchanges $A_{u,1}$ and $A_{u,2}$ as well as $E_1$ and
  $E_2$. Hence $\iota_{\omega}$ interchanges $A_{u,1}$ and $A_{u,2}$ but
  induces non-trivial involutions on $E_1$ and $E_2$. In order to
  determine the fixed points of $\iota_{\omega}$ it is sufficient to
  compute the limit of the fixed points of $\iota_{\omega}$ in $A_{\tau}$
  as $\tau_1\rightarrow i\infty$. The sixteen fixed points of $\iota_{\omega}$
  on $A_{\tau}$ are given by
$$
\left[(\tau_1+\tau_2,
\tau_2+\tau_3)/2+\varepsilon_1(\tau_1,\tau_2)+\varepsilon_2(\tau_2,\tau_3)+
\varepsilon_3
(1,0)+ \varepsilon_4(0,3)\right]\in A_{\tau}
$$
where $\varepsilon_i=0$ or $1$. As $\tau_1\rightarrow i\infty$ these 16
points come together in pairs; more precisely any two points which only
differ by $\varepsilon_3$ have the same limit. This gives us eight points of
which four lie on each of the curves $E_i$ (depending on whether
$\varepsilon_1=0$ or $1$). These points are given by

$$
\begin{array}{ll}
\left[(\tau_2+\tau_3)/2+\varepsilon_2\tau_3+\varepsilon_4 3\right] \in
E(\tau_3) &  (\varepsilon_1=0)\\[2mm]
\left[(\tau_2+\tau_3)/2+\tau_2+\varepsilon_2\tau_3+\varepsilon_4
3\right]\in E(\tau_3) & (\varepsilon_1=1).
\end{array}
$$
\hfill
\end{proof}

Figure~$2$ indicates the action of $\iota_{\omega}$ and the position of the
eight fixed points on $\bar A_u$.

\unitlength1cm
\begin{figure}[htb]
\begin{center}
{\includegraphics[width=0.5\columnwidth]{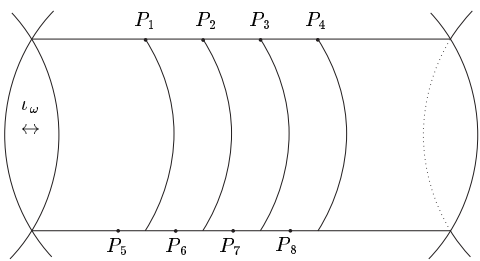}}
\end{center}
\caption{The involution $\iota_{\omega}$}
\end{figure}

Note that two fixed points lie on one ruling if and only if $[\tau_2]$ is a
$2$-torsion point on $E(\tau_3)$, i.e. if and only if the glueing parameter
$e=[2\tau_2]=0$.

The next step is to extend the polarization to the degenerate abelian
surfaces. Ideally we would like to glue the line bundle ${\cal L}_{\tau}$
on $A_{\tau}$ to a line bundle ${\cal L}$ on $A$ and to extend this line
bundle to ${\bar A}$ in such a way that the sections ${\widehat s}_{\alpha
\beta}$ as well as the action of the symmetry group (see Proposition \ref
{prop12}) extend. At this point, however, we encounter a fundamental
difficulty. We have seen that it is possible to extend $A$ to ${\bar A}$ in
such a way that the symmetries, and here in particular the involution
$\iota_{\omega}$, extend. It is also possible to define a suitable line
bundle ${\cal L}$ and extend it to a line bundle ${\bar {\cal L}}$ on
${\bar A}$ (see \cite[II.5]{HKW}). But it is not possible to do this in
such a way that the action of the symmetry group also extends to
${\bar{\cal L}}$. (This leads in particular to a numerical contradiction on
the fibre over the origin $0\in B$.) For this reason we shall now restrict
ourselves to taking the partial quotient with respect to the group
$$
P'=\left\{
\left(
\begin{array}{c|c}
{\bf 1}   &
\begin{array}{cc}
2\ZZ\ & 0\\
0     & 0
\end{array}\\
\hline
0         &   {\bf 1}
\end{array}
\right)
\right\}
\subset \mbox {Sp}(4,\ZZ).
$$
This is the lattice contained in the stabilizer of the line generated by
$l_0=(0, 0, 1, 0)$ in the group $\Gamma_{1,3}(2)\cap
\Gamma_{1,3}^{\mbox{\scriptsize lev}} $. The partial quotient defined by
this group is given by the map
$$
\begin{array}{rcl}
\HH_2 & \rightarrow &
\CC^*\times\CC\times\HH_1\subset\CC\times\CC\times\HH_1\\[2mm]
\left(\begin{array}{cc}
\tau_1  &  \tau_2\\
\tau_2  &  \tau_3
\end{array}\right)
& \mapsto & (t_1=e^{2\pi i\tau_1/2}, \tau_2, \tau_3).
\end{array}
$$
Partial compactification of $\HH_2/P'$ in $\CC\times\CC\times\HH_1$ is given by
$$
B':=\overset{\circ}{(\overline{\HH_2/P'})}
$$
The two partial quotients with respect to $P$ and $P'$ are related by
the glueing map
$$
\begin{array}{rcl}
\varphi:\quad B'  &  \rightarrow  &  B\\
(t_1, \tau_2, \tau_3)  &  \mapsto & (t_1 t_2, t_2 t_3, t_2^{-1})
\end{array}
$$
where $t_2=e^{2\pi i\tau_2/6}$ and $t_3=e^{2\pi i \tau_3/18}$. The image
of $\varphi$ is $B\backslash (B\cap\{T_2T_3=0\})$ and the map $\varphi$ is
unramified onto its image. We can pull the family ${\bar A}$ over $B$ back
to $B'$ via $\varphi$ and we shall denote the resulting family by ${\bar
A'}$. This family extends the universal family $A'$ over $\HH_2/P'$. We
shall denote the projection from ${\bar A'}$ to $B'$ by $\pi'$.

\begin{proposition}\label{prop23}
  {\rm(i)} The line bundles ${\cal L}_{\tau}$ on $A_{\tau}$ glue to a line
  bundle ${\cal L}'$ on~$A'$.\\
  {\rm(ii)} The line bundle ${\cal L}'$ can be extended to a line bundle
  ${\bar{\cal L'}}$ on ${\bar A'}$ in such a way that the sections
  ${\widehat s}_{\alpha \beta}$ as well as the action of the Heisenberg
  group $H_{2 6}$ and the involution $\iota_{\omega}$ extend.
\end{proposition}

\begin{proof}
  A straightforward computation shows that with respect to the coordinates
  $t_1=e^{2\pi i\tau_1/2}$ and $w_1=e^{2\pi i z_1/2}, w_2=e^{2\pi i z_2/6}$
  the functions ${\widehat\Theta}_{\alpha \beta}(\tau, z)$ are given by:
$$
\begin{array}{rcl}
{\widehat\Theta}_{\alpha \beta}(\tau, z) & = &
{\Theta}_{00}\frac{\alpha}{2}\frac{\beta}{6}(\tau', z' -{\omega}')\\[2mm]
& = & \sum\limits_{q\in\ZZ^2} t_1^{\frac 12 q_1(q_1-1)} \exp\{9\pi i
q_2(q_2-3)\tau_3\}\\[3mm]
& & \qquad\cdot\exp\{6\pi i(2q_1 q_2 -3 q_1 -q_2)\tau_2\}w_1^{q_1}
w_2^{q_2} (-1)^{\alpha q_1}\rho_6^{\beta q_2}.
\end{array}
$$
In particular this shows that we can consider these functions as
functions on $(\HH_2/P')\times (\CC^*)^2$. Similarly we find that with
respect to the lattice $L_{\tau}$ the functions ${\widehat\Theta}_{\alpha
\beta}(\tau, z)$ have the following transformation behaviour. For $(k,
l)\in \ZZ^2$:
$$
{\widehat\Theta}_{\alpha\beta} (\tau,
z+(2k,6l))={\widehat\Theta}_{\alpha\beta}(\tau,z).
$$
For $(m,n)\in \ZZ^2$:
$$
\begin{array}{l}
{\widehat\Theta}_{\alpha\beta} \left(\tau, z+(m,n)\left(
\begin{array}{cc}
2\tau_1 & 2\tau_2\\
2\tau_2 & 2\tau_3
\end{array}
\right)\right)  \\
\qquad= \Theta_{00\frac{\alpha}{2}\frac{\beta}{6}}\left(\tau',
z'-{\omega}'+(2m,6n)
\left(
\begin{array}{cc}
\tau'_1 & \tau'_2\\
\tau'_2 & \tau'_3
\end{array}\right)\right)\\[2mm]
\qquad=\exp\left\{2\pi i\left[-\frac 12 (2m, 6n)\left(
\begin{array}{cc}
\tau'_1 & \tau'_2\\
\tau'_2 & \tau'_3
\end{array}\right)\left(\begin{array}{c} 2m\\ 6n\end{array}\right)-(2m,
6n)(z'-{\omega}')\right]\right\}\\[2mm]
\qquad\qquad\cdot\Theta_{2m, 6n,\frac {\alpha}{2},\frac{\beta}{6}} (\tau',
z'-{\omega}')\\[2mm]
\qquad= t_1^{-2m^2+m} w_1^{-2m} w_2^{-6n} e^{2\pi
i[(-2mn+\frac{m^2}{2})\tau_2+(-n^2+\frac n2)\tau_3]}
{\widehat\Theta}_{\alpha \beta}(\tau, z).
\end{array}
$$
\noindent
These calculations show claim (i). To prove (ii) we have to consider the
limit as $t_1\rightarrow 0$. Here we find
$$
\begin{array}{rcl}
\lim\limits_{t_1\rightarrow 0}{\widehat\Theta}_{\alpha \beta}(\tau,
z) & = &\sum\limits_{q_2\in \ZZ} e^{9\pi i q_2(q_2-3)} e^{6\pi i
(-q_2\tau_2)}
w_2^{q_2} \rho_6^{\beta q_2}\\
&&\quad+(-1)^{\alpha} w_1\sum\limits_{q_2\in\ZZ} e^{9\pi i q_2(q_2-3)}
e^{2\pi i
(q_2-3)\tau_2} w_2^{q_2} \rho_6^{\beta q_2}\\
& = &\vartheta_{0\frac k6} \left(\tau_3 /
6,(z_2-\tau_3/2-\tau_2/2)/6\right )
\\
&&\quad+(-1)^{\alpha} w_1 e^{2\pi
i\left(-\frac{\tau_2}{4}\right)}\vartheta_{0\frac k6}\left(\tau_3
/6,(z_2-\tau_3/2+\tau_2/2)/6
\right).
\end{array}
$$
To prove that ${\cal L}'$ can be extended to a line bundle ${\bar{\cal
L}'}$ on $\bar A'$ we can argue exactly as in the proof of
\cite[Proposition (II.5.13)]{HKW}, the only difference being that we took
the partial quotient with respect to a smaller group. (Note that if we take
the quotient with respect to $P$ we no longer obtain integer exponents of
$t_2$.) The extension of the action of the symmetry group follows as in the
proof of \cite[Proposition (II.5.41)]{HKW}.  \hfill
\end{proof}

\section {The map to $\PP^3$}

We consider boundary points $u'=(0,\tau_2, \tau_3) \in B'$ and
$u=(0,t_2t_3, t_2^{-1})\in B$ (so $u=\varphi(u')$) and the associated
degenerate abelian surfaces
$$
{\bar A'}_{u'}={\bar A}_u=A_{u,1}\cup A_{u,2}
$$
where $A_{u,1}=A_{u,2}$ is an elliptic ruled surface. We gave a precise
description of the surfaces $A_{u,i}$ and the way they are glued in Theorem
\ref{theo21}. Recall that $A_{u,i}$ is an elliptic ruled surface with two
disjoint sections $E_i, i=1,2$ of self-intersection number $E_i^2=0$. We
shall denote the fibre over a point $P$ of the base curve by $f_P$. Recall
also the line bundle ${\bar{\cal L}}'$ on ${\bar A}'$. We set
$$
{\bar{\cal L}}_u:={\bar{\cal L}}'|_{{\bar A}'_u}={\bar{\cal L}}'|_{{\bar A}_u},
\quad {\cal
L}_{u,i}:={\bar{\cal L}}'|_{A_{u,i}}.
$$
\noindent
In the proof of Proposition \ref{prop23} we computed that
$$
\begin{array}{rl}
\lim\limits_{t_1\rightarrow 0}
{\widehat\Theta}_{\alpha\beta}(\tau,z)=&
\vartheta_{0\frac{k}{6}}(\tau_3/6,(z_2-\tau_3/2-\tau_2/2)/6)\\[2mm]
&+w^{-1}\vartheta_{0\frac{k}{6}}(\tau_3/6,(z_2-\tau_3/2+\tau_2/2)/6).
\end{array}
$$
The theta function $\vartheta_{0\frac k6}(\tau_3/6,
(z_2-\tau_3/2+\tau_2/2)/6)$ has six zeroes on the elliptic curve
$E(\tau_3)=\CC/(\ZZ 2\tau_3+\ZZ6)$. Hence
$$
\operatorname{deg} {\cal L}_{u,i}|_{E_i}=6.
$$
Since the exponent of $w$ is $-1$ it follows that
$$
\operatorname{deg} {\cal L}_{u,i}|_{f_P}=1.
$$
(See \cite[Proposition(II.5.35)]{HKW} for similar considerations in
the $(1,p)$ case.) Hence
$$
{\cal L}_{u,i}={\cal O}_{A_{u,i}}(E_1+6 f_{P_i})
$$
for a suitable point $P_i\in E(\tau_3)$. (The point $P_i$ can be
computed from $\lim\limits_{t_{1\rightarrow
0}}{\widehat\Theta}_{\alpha\beta}(t,z)$ and the normal bundle of $E_1$ in
$A_{u,i}$, but we shall not need this later.) Standard arguments using
Riemann-Roch show that
$$
h^0(A_{u,i}, {\cal L}_{u,i})=h^0(A_{u,i}, {\cal O}_{A_{u,i}}(E_1+6f_{P_i}))=12.
$$
\begin{proposition}\label{prop31}
  The restriction ${\bar{\cal L}}_u={\bar{\cal L}}'|_{{\bar A}'_{u'}}
  ={\bar{\cal L}}'|_{{\bar A}_{u}}$ of ${\bar{\cal L}}'$ to the degenerate
  abelian surface ${\bar A}_u$ has the
  following properties:\\
  {\rm (i)} $h^0({\bar A}_u,{\bar{\cal L}}_u)=12,$\\
  {\rm (ii)} The restriction map $\operatorname{rest}: H^0({\bar
  A}',{\bar{\cal
  L}}')\rightarrow H^0({\bar A}_u, {\bar{\cal L}}_u)$ is surjective.\\
  {\rm (iii)} The restriction map $\operatorname{rest}: H^0({\bar
  A}_u,{\bar{\cal L}}_u)\rightarrow H^0 (A_{u,i} {\cal L}_{u,i})$ is an
  isomorphism.
 \end{proposition}

\begin{proof}
  We consider the space $V\subset H^0({\bar A}',{\bar{\cal L}}')$ which is
  spanned by the twelve sections ${\widehat s}_{\alpha \beta}$. Since these
  sections are a basis of $H^0(A_{\tau},{\cal L}_{\tau})$ for every $\tau$
  the space $V$ has dimension~$12$. We claim that the restriction map
$$
\mbox{rest}:  V\rightarrow H^0({\bar A}_u, {\bar{\cal L}}_u)
$$
is injective. By our computation of $\lim\limits_{t\rightarrow 0}
{\widehat\Theta}_{\alpha \beta}(\tau, z)$ it follows that this map is not
identically zero. It is also $H_{26}$-equivariant and hence our claim
follows if we can show that $V$ is irreducible as an $H_{26}$--module. But
this is easy to see: As an $H_6$-module $V=V_0\oplus V_1$ where
$V_i=\mbox{span } (\widehat s_{i\beta}, \beta=0,\ldots, 5)$. The
$H_6$--modules $V_0$ and $V_1$ are irreducible, Moreover addition by
$e_3/2$ interchanges $V_0$ and $V_1$. Hence $h^0({\bar A}_u, {\bar{\cal
L}}_u)\ge 12$.

We have already remarked that $h^0(A_{u,i},{\cal
O}_{A_{u,i}}(E_1+6f_{P_i}))=12$. Our next claim is that the map
$$
\mbox{rest }: H^0(A_{u,i}, {\cal L}_{u,i})\rightarrow H^0(E_1,{\cal
L}_{u,i}|_{E_1})\oplus H^0(E_2,{\cal L}_{u,i}|_{E_2})
$$
is an isomorphism. Since the vector spaces on both side have the same
dimension, namely~$12$, it is enough to prove injectivity. This follows
from $h^0(A_{u,i},{\cal O}_{A_{u,i}}(-E_2+6f_{P_i}))=0$. But now this
implies that glueing sections on $A_{u,1}$ and $A_{u,2}$ along $E_1$ and
$E_2$ gives at least $2\times 6=12$ conditions. Hence $h^0({\bar A}_u,
{\bar{\cal L}}_u)\le 12$. With our previous argument this shows that
$h^0({\bar A}_u, {\bar{\cal L}}_u)=12$ and hence both (i) and (ii) are
proved. This also shows that
$$
\mbox{rest}: H^0({\bar A}_u,{\bar{\cal L}}_u)\rightarrow H^0(A_{u,i}, {\cal
L}_{u,i})
$$
is an isomorphism and hence we have shown (iii).\hfill
\end{proof}
Since we are interested in the map to $\PP^3$ given by $H^0({\cal L})^-$ we
consider the subspace
$$
V^-=\langle {\widehat g}_0, {\widehat g}_1, {\widehat g}_2, {\widehat
g}_3\rangle\subset V
\subset H^0({\bar A}', {\bar{\cal L}}')
$$
and
$$
\begin{array}{lcl}
V^-_{u}&=&\mbox{ rest } (V^-\rightarrow H^0 ({\bar A}_u, {\bar {\cal
L}}_u)),\\[2mm]
 V^-_{u,i}&=&\mbox{ rest } (V^-\rightarrow H^0 (A_{u,i},
{\cal L}_{u,i})).
\end{array}
$$
The spaces $V^-_u$ and $V^-_{u,i}$ are $4$--dimensional. We want to
study the map
$$
\varphi_{V^-_u}:{\bar A}_u--\rightarrow \PP^3.
$$
The sections ${\widehat g}_i$ vanish at the eight ``$2$-torsion'' points
$P_1,\ldots ,P_8$ on ${\bar A}_u$ and hence
$$
V^-_u\subset H^0(A_{u,i}, {\cal O}_{A_{u,i}}(E_1+6f_{P_i}-\sum\limits^8_{j=1}
P_j)).
$$
Again using restriction to $E_1$ and $E_2$ it follows that the vector
space on the right hand side has dimension~$4$ and hence
$$
V^-_u=H^0(A_{u,i},{\cal O}_{A_{u,i}}( E_1+6f_{P_i}-\sum\limits^8_{j=1} P_j)).
$$
We shall first consider the {\em product case}, i.e. $e=[\tau_2]=0$. In
this case $A_{u,i}=E(\tau_3)\times \PP^1$ and there are four rulings which
contain two of the points $P_j$ each. These four rulings are, therefore, in
the base locus of the linear system $|V^-_u|$. Removing this base locus we
obtain the complete linear system of a line bundle on $A_{u,i}$ which has
degree~$2$ on the sections $E_i$ and degree $1$ on the fibres. This maps
$A_{u,i}$ $2:1$ onto a quadric. Since the map ${\bar A}_u--\rightarrow
\PP^3$ factors through $A_{u,i}$ this shows that the ``Kummer surface''
${\bar A}_u /\iota_{\omega}$ is mapped $2:1$ onto a quadric. It should be
noted that double quadrics arise not only from degenerations of abelian
surfaces, but also from special abelian surfaces, namely products (cf.
Theorem \ref{theo13}). In fact what happens is that the map from the moduli
space ${\cal A}_{1,3}(2)$ (or its extension to a toroidal
compactification) contracts each Humbert surface parametrizing product
surfaces to a double point in $N$ corresponding to a quadric.

{}From now on we shall assume $e\neq 0$. Then no fibre of the ruling
contains two of the points $P_i$. We have to recall the notion of an {\em
elementary transformation} of a ruled surface $S$ at a point $P$. This
consists of first blowing up $S$ in $P$ and then blowing down the strict
transform of the fibre through $P$. The result is again a ruled surface
$\mbox{elm}_{P}S$. Let
$$
{\hat A}_u:=\mbox{elm}_{P_1,\ldots,P_8}(A_{u,i}).
$$
Then ${\hat A}_u$ has again two disjoint sections $E_1$ and $E_2$ with
$E_i^2=0$. In particular
$$
{\hat A}_u=\PP({\cal O}_{E(\tau_3)}\oplus{\widehat{\cal M}}_u) \mbox { for
some }
{\widehat{\cal M}}_u\in \mbox{Pic}^{0}(E(\tau_3)).
$$
Note that ${\widehat{\cal M}}_u$ and ${{\widehat{\cal M}}_u}^{-1}$
define the same $\PP_1$-bundle. It is straightforward to compute the normal
bundle of the sections $E_1$ and $E_2$ in $A_{u,i}$ (cf. \cite[p. 229]
{HKW}. Since we can control the self-intersection of a section under
blowing up and blowing down and since we know the points $P_j$ it is
straightforward to show that ${{\widehat{\cal M}}_u}^{\pm 1}={\cal
O}_{E(\tau_3)}(2[\tau_2]-2[0])$, and hence
$$
{\hat A}_u=\PP({\cal O}_{E(\tau_3)}\oplus {\cal
O}_{E(\tau_3)}(2[\tau_2]-2[0])).
$$
Consider the diagram
$$
\unitlength1pt
\begin{picture}(100,60)(0,0)
\put(45,55){$\tilde A_u$}
\put(16,35){$\pi_1$}
\put(75,35){$\pi_2$}
\put(43,51){\vector(-1,-1){35}}
\put(53,51){\vector(1,-1){35}}
\put(0,5){$A_u$}
\put(90,5){$\hat A_u$}
\end{picture}
$$
where $A_u=A_{u,i}$ and ${\tilde A}_u$ is $A_u$ blown up in $P_1,\ldots,
P_8$. We denote the exceptional divisors over $P_1,\ldots, P_8$ by
$E^1,\ldots, E^8$. The line bundle $\pi^*_1{\cal L}_{u,i}\otimes{\cal
O}_{{\tilde A}_u}(-E^1-\ldots-E^8)$ has degree $0$ on the strict transforms
of the fibres through the points $P_j$. Hence
$$
{\hat{\cal L}}_u:=\pi_{2*}(\pi^*_1{\cal L}_{u,i}\otimes {\cal O}_{{\tilde
A}_u}(-E^1-\ldots-E^8))\in \mbox{Pic } {\hat A}_u
$$
is a line bundle on ${\hat A}_u$. The degree of ${\hat{\cal L}}_u$ is
$1$ on a ruling and $2$ on the sections $E_i$. Hence
$$
{\hat{\cal L}}_u={\cal O}_{{\hat A}_u}(E_1+2f_Q)
$$
for a suitable point $Q$ on the base curve. (Clearly $Q$ can be computed
explicitly, but this is immaterial for our purposes.) By the usual
arguments
$$
h^{0}({\hat A}_u, {\hat{\cal L}}_u)=4
$$
and the rational map from $A_u$ to ${\hat A}_u$ defines an isormorphism
$$
\pi_{2*} \pi^*_1:V^-\cong H^{0}({\hat A}_u, {\hat{\cal L}}_u).
$$
\begin{proposition}\label{prop32}
  Let $e\neq 0$. The linear system $|V^-|$ on ${\hat A}_u$ has the
  following
  properties:\\
  {\rm(i)} $|V^-|$ is base point free.\\
  {\rm(ii)} $|V^-|$ maps the two sections $E_i$ each $2:1$ onto two skew
  lines.\\
  {\rm(iii)} $|V^-|$ is very ample outside the sections $E_i$. More
  precisely, if a cluster $\zeta$ of length~$2$ (i.e. two points or a point
  and a tangent direction) is not embedded, then $\zeta$ is contained in
  $E_1$ or $E_2$.
\end{proposition}

\begin{proof}
This follows easily from Reider's theorem. We write
$$
|V^-|=|E_1+2f_Q|=|K+L|
$$
where the canonical divisor $K=-E_1-E_2$ and $L=2E_1+E_2+2f_Q$. Then
$L^2=12$. If $|V^-|$ is not base point free, then there exists a curve $D$
with $L.D=0$ or $1$. Clearly such a curve cannot exist. If $|V ^-|$ fails
to embed a cluster $\zeta$ then there exists a curve $D\supset\zeta$ with
$L.D=2$ and $D^2=0$. Then $D$ must be a section with $D.E_i=0$. Since
${\widehat{\cal M}}_u\neq{\cal O}_{E(\tau_3)}$ (here we use
$e=2[\tau_2]\neq 0)$ it follows that $D=E_1$ or $D=E_2$. Finally note that
the restriction of $|V^-|$ to the elliptic curves $E_i$ gives a complete
linear system of degree~$2$. Hence these curves are mapped $2:1$ to lines.
Since pairs $(x,y)$ with $x\in E_1$ and $y\in E_2$ are separated, there
lines are skew.  \hfill
\end{proof}
We can now summarize our results as follows.
\begin{theorem}\label{theo33}
  Let ${\bar A}_u$ be a corank~$1$ degenerate abelian surface over a point
  $u=(0,T_2,T_3)\in B$ with $T_2T_3\neq0$ and consider the Kummer map given
  by the linear system $|V^-|$:
$$
\phi_{|V^-|}:{\bar A}_u--\rightarrow\PP^3.
$$
{\rm(i)} If $e=0$ then ${\bar A}_u$ is mapped $4:1$ onto a smooth quadric.\\
{\rm(ii)} If $e\neq 0$ then there is a commutative diagram.
$$
\unitlength1pt
\begin{picture}(100,60)(0,0)
\put(0,46){$\bar A_u$}
\put(13,50){\line(1,0){15}}
\put(43,50){\line(1,0){15}}
\put(73,50){\vector(1,0){20}}
\put(100,46){$\PP^3$}
\put(13,45){\line(1,-1){10}}
\put(27,31){\line(1,-1){10}}
\put(12,20){$2:1$}
\put(41,17){\vector(1,-1){12}}
\put(54,-3){$\hat A_u$}
\put(67,8){\vector(1,1){30}}
\end{picture}
$$
The image of ${\bar A}_u$ is an elliptic ruled surface which has double
points along two skew lines but no other singularities.
\end{theorem}

\begin{proof}
  By construction the Kummer map ${\bar A}_u--\rightarrow \PP^3$ factors
  through ${\bar A}_u/{\iota_{\omega}}=A_{u,i}$. All other statements
  follow from our above discussion of the map $\phi_{|V^-|}: {\hat A}_u \to
  \PP^3$ given by the linear system $|V^-|$.  \hfill
\end{proof}

This theorem explains the irreducible singular quartic surfaces which are
parametrized by the S-planes, appearing already in \cite[Section 5-2]{Ni}.

We want to conclude this paper with some remarks. We have already seen that
the degenerate surfaces with $e=0$ correspond to products. The limits of
{\em bielliptic} abelian surfaces are characterised by $2e=0, e\neq 0$.
Geometrically this means that the surfaces ${\bar A}_u$ contain degenerate
elliptic curves which are $4$-gons, i.e. cycles consisting of $4$ rulings.
The degenerate abelian surfaces ${\bar A}_u$ where $u$ is on another
boundary component or where more than one of the $T_i$ vanishes can also be
described. If $T_3=0$ and $T_1 T_2\neq 0$ then ${\bar A}_u$ is a chain
of~$6$ elliptic ruled surfaces. If two of the $T_i$ are zero, then ${\bar
A}_u$ consists of $12$ quadrics, whereas ${\bar A}_0$ has $36$ components,
of which $24$ are $\PP^2$ and $12$ are $\PP^2$ blown up in three points.
Limits of polarizations of type $(2,6)$ exist on these surfaces, but as we
pointed out before, there is no possibility of defining the Kummer map
globally over $B$. On the other hand it is easy to construct degenerations
of the ruled surfaces $A_{u,i}$ which lead to a union of two quadrics
intersecting along a quadrangle or to a tetrahedron.

Finally we want to comment on the boundary components of the Igusa
compactification of the moduli space ${\cal A}_{1,3}(2)$. These are
enumerated by the {\em Tits building} of the group $\Gamma_{1,3}(2)$, i.e.
by the equivalence classes modulo $\Gamma_{1,3}(2)$ of the lines and
isotropic planes in $\QQ^4$. The Tits building was calculated by Friedland
in \cite{F}: for details, and for some other cases, see \cite{FS}.
There are $30$ equivalence classes of lines. These correspond to the $15$
equivalence classes of short, respectively long vectors. Each set of $15$
lines is naturally parametrized by $(\ZZ/2)^4 \setminus
\{0\}=\PP^3(\FF_2)$. The $15$ planes are parametrized by the $15$ isotropic
planes in $\operatorname{Gr}(1,\PP^3(\FF_2))$. The isotropic planes are a
hyperplane section of $\operatorname{Gr}(1,\PP^3(\FF_2))$ embedded as a
quadric via the Pl\"ucker embedding. The $15$ short and the $15$ long
vectors as well as the $15$ planes are identified under the group
$\Gamma_{1,3}/\Gamma_{1,3}(2)\cong \operatorname{Sp}(4,\FF_2)\cong S_6$.
That is, there are two equivalence classes of lines modulo $\Gamma_{1,3}$
and one plane (see also \cite [Theorem(I.3.40)]{HKW}). Finally the
involution $V_3$ of the maximal arithmetic subgroup $\Gamma^*_{1,3}$
identifies short and long vectors (see \cite[Folgerung 3.7]{G} and
\cite[Section 2]{HNS}). In our computations above the boundary component
given by $T_1=0$ corresponds to a short vector, whereas the boundary
component given by $T_3=0$ corresponds to a long vector. We described the
degenerate abelian surfaces associated to points on a boundary component
correspronding to a short vector. The matrix $V_3$ (and similarly any
involution $gV_3$ where $g$ is an element of $\Gamma_{1,3}$ -- cf.
\cite[Theorem 2.4]{HNS}) interchanges boundary components associated to
short vectors with boundary components associated to long vectors. It
should, however, be pointed out that the induced action of $V_3$ on the
Igusa compactification ${\cal A}_{1,3}^*$ is only a rational map, not a
morphism. This follows since the boundary components associated to long,
and short vectors are not isomorphic: although their open parts (i.e.
away from the corank-$2$ boundary components) are isomorphic (namely to the
open Kummer modular surface $K^{0}(1)$), they contain different
configurations of rational curves in the corank-$2$ boundary components.
This follows from \cite[Satz III.5.19]{B} and \cite[Theorem 4.13]{W}. The
degenerate abelian surfaces belonging to points on $T_3=0$ are different
from those associated to points on $T_1=0$: they are a cycle of six
elliptic ruled surfaces rather than two. At first this looks like a
contradiction to \cite[Theorem 2.4]{HNS}, but this is not the case. The
polarization on each of the six components of the surface $A_P$ where $P\in
\{T_3=0\}$ is of the form ${\cal O}(E_1+2f_P)$. On four of the six
components the linear system $|V^-|$ has a base locus consisting of a
section. These four components are contracted. The other two components are
identified under $|V^-|$ and are mapped to a quartic which is an elliptic
ruled surface singular along two skew lines. In this way we find the same
images in $\PP^3$ as in the case $T_1=0$.

%
%
\bibliographystyle{amsalpha}

\end{document}